\title{Regular dessins with a given automorphism group}
\author{Gareth A. Jones\\
School of Mathematics\\
University of Southampton\\
Southampton SO17  1BJ, U.K.\\
{\tt G.A.Jones@maths.soton.ac.uk}
}

\documentclass[10pt]{article}
\usepackage{color}
\usepackage[latin1]{inputenc}
\usepackage{a4wide}
\usepackage{latexsym}
\usepackage{amsmath}
\usepackage{amssymb}
\usepackage{amsfonts}
\newtheorem{thm}{Theorem}[section]
\newtheorem{lemma}[thm]{Lemma}
\newtheorem{cor}[thm]{Corollary}

\date{}
\begin{document}

\maketitle


\medskip

\begin{abstract}
Dessins d'enfants are combinatorial structures on compact Riemann surfaces defined over algebraic number fields, and regular dessins are the most symmetric of them. If $G$ is a finite group, there are only finitely many regular dessins with automorphism group $G$. It is shown how to enumerate them, how to represent them all as quotients of a single regular dessin $U(G)$, and how certain hypermap operations act on them. For example, if $G$ is a cyclic group of order $n$ then $U(G)$ is a map on the Fermat curve of degree $n$ and genus $(n-1)(n-2)/2$. On the other hand, if $G=A_5$ then $U(G)$ has genus $274218830047232000000000000000001$. For other non-abelian finite simple groups, the genus is much larger.
\end{abstract}

\medskip

 (secondary).

\section{Introduction}

Bely\u\i's Theorem~\cite{Bel} states that a compact Riemann surface $C$, regarded as a complex projective algebraic curve, can be defined over the field $\overline{\mathbb Q}$ of algebraic numbers if and only if there is a non-constant meromorphic function $\beta:C\to{\mathbb P}^1({\mathbb C})$ ramified over at most three points. This is equivalent to the complex structure on $C$ being obtained, in a canonical way, from a combinatorial structure $D$ on $C$ called a dessin; this can be regarded as an oriented hypermap, a bipartite map, or a tripartite triangulation of $C$. The most symmetric dessins are the regular dessins, those for which $\beta$ is a regular covering, in which case the automorphism group ${\rm Aut}\,D$ of $D$, the group of covering transformations of $\beta$, and the monodromy group of $\beta$ are all isomorphic to a single $2$-generator finite group $G$. In principle, understanding regular dessins is sufficient for an understanding of all dessins, since each dessin arises as the quotient of some regular dessin by a group of automorphisms.

The set ${\mathcal R}(G)$ of regular dessins associated in this way with a given $2$-generator  finite group $G$ corresponds bijectively to the set ${\mathcal N}(G)$ of normal subgroups $N$ of a free group $F$ of rank $2$ with $F/N\cong G$, or equivalently to the set of orbits of ${\rm Aut}\,G$ on generating pairs of $G$. Our aim here is to show how a method due to Hall~\cite{Hal} can be used to determine the number $r(G)=|{\mathcal R}(G)|$ of such dessins, to show how they can be represented as quotients of a single regular dessin $U(G)$, and to describe the automorphism group $\overline G={\rm Aut}\,U(G)$ of $U(G)$. A number of illustrative examples are given, including cases where $G$ is a cyclic group, a dihedral group or a non-abelian finite simple group; in the latter case, $U(G)$ and $\overline G$ turn out to be surprisingly large.

The absolute Galois group $\Gamma={\rm Gal}\,\overline{\mathbb Q}/{\mathbb Q}$ acts on dessins by acting on the coefficients of the polynomials and rational functions defining them over $\overline{\mathbb Q}$. By a recent result of Gonz\'alez-Diez and Jaikin-Zapirain~\cite{GJ}, it acts faithfully on regular dessins. It preserves automorphism groups, so it acts on ${\mathcal R}(G)$ for each $G$. A group $\Omega\cong {\rm Out}\,F\cong GL_2({\mathbb Z})$ of hypermap operations, introduced by James~\cite{Jam}, also preserves regularity and automorphism groups, and its orbits on ${\mathcal R}(G)$ correspond to the $T$-systems of generating pairs for $G$. We shall consider some examples of such actions of $\Gamma$ and $\Omega$, interpreting group-theoretic results of Neumann and Neumann~\cite{NN}, Dunwoody~\cite{DunTS}, and Garion and Shalev~\cite{GS} on $T$-systems in terms of dessins.

\medskip

\noindent{\bf Acknowledgements} The author is very grateful to Martin Dunwoody, Pierre Guillot and Frank Herrlich for their helpful comments about $T$-systems, the absolute Galois group and origamis, and to the organisers of the conference on Riemann and Klein Surfaces, Symmetries and Moduli Spaces, in honour of Emilio Bujalance (Link\"oping, June, 2013), for the invitation to give a talk on which this paper is based.

\section{Dessins}\label{des}

If $C$ is a compact Riemann surface, or equivalently a nonsingular complex projective algebraic curve, then by Bely\u\i's Theorem~\cite{Bel}, as reinterpreted and extended by Grothendieck~\cite{Gro}, Wolfart~\cite{Wol} and others, the following are equivalent:

\begin{enumerate}
\item $C$ is defined, as an algebraic curve, over the field $\overline{\mathbb Q}$;
\item there is a non-constant meromorphic function $\beta:C\to{\mathbb P}^1({\mathbb C})$, branched over at most three points;
\item $C$ is uniformised by a subgroup $M$ of finite index in a triangle group $\Delta$;
\item $C$ is isomorphic to  the compactification $\overline{{\mathbb H}/N}$, where $\mathbb H$ is the hyperbolic plane and $N$ is a subgroup of finite index in the congruence subgroup $\Gamma(2)$ of level $2$ in the modular group $PSL_2({\mathbb Z})$.
\end{enumerate}
The function $\beta$ in~$(2)$, called a {\it Bely\u\i\/ function}, corresponds to the surface coverings induced by the inclusions $M\le\Delta$ and $N\le\Gamma(2)$ in $(3)$ and $(4)$; by applying a M\"obius transformation, the ramification points of $\beta$ can, without loss of generality, be taken to be contained in $\{0, 1, \infty\}$. The group $\Gamma(2)$ is freely generated by the parabolic M\"obius transformations
\[X:z\mapsto\frac{z}{-2z+1}\qquad{\rm and}\qquad Y:z\mapsto\frac{z-2}{2z-3}\]
fixing $0$ and $1$, and there is an epimorphism $\Gamma(2)\to\Delta$, with $N$ the inverse image of $M$. The compactification in~$(4)$ is obtained by adding points to fill the punctures in ${\mathbb H}/N$ corresponding to the cusps of $N$, that is, its orbits on ${\mathbb P}^1({\mathbb Q})$. Further background is given by Girondo and Gonz\'alez-Diez in~\cite{GG}, and by Lando and Zvonkin in~\cite{LZ}.

Under the above conditions, $C$ is called a {\it Bely\u\i\/ curve}. It carries a {\it dessin}, that is, a bipartite map in which the black and white vertices are the fibres of $\beta$ over $0$ and $1$, and the edges are the inverse images of the unit interval $I=[0,1]$. The {\it monodromy group} $G$ of $D$ is the monodromy group of $\beta$, restricted to a smooth covering of ${\mathbb P}^1({\mathbb C})\setminus\{0, 1, \infty\}$; this is the transitive group of permutations of the sheets of this covering obtained by lifting closed paths in ${\mathbb P}^1({\mathbb C})\setminus\{0, 1, \infty\}$ to $C$. Equivalently, $G$ is the group of permutations of the edges of $D$ generated by the monodromy permutations $x$ and $y$ about $0$ and $1$, obtained by following the orientation of $C$ around the black and white vertices; this is isomorphic to the group induced by $\Delta$ on the cosets of $M$, or equivalently by $\Gamma(2)$ on the cosets of $N$, with $X$ and $Y$ inducing $x$ and $y$.

The automorphism group ${\rm Aut}\,D$ of a dessin $D$ (preserving orientation and vertex colours) is the group of permutations of the edges commuting with $x$ and $y$, or equivalently with the monodromy group $G$. It can be identified with the group of covering transformations of $\beta$, permuting the fibre over a base-point. It acts semiregularly on the edges, and we say that $D$ is a {\it regular dessin} if it acts transitively (and hence regularly). This is equivalent to $M$ and $N$ being normal subgroups of $\Delta$ and $\Gamma(2)$, in which case
\[{\rm Aut}\,D\cong G\cong\Delta/M\cong\Gamma(2)/N,\]
with ${\rm Aut}\,D$ and $G$ acting on edges as the left and right regular representations of the same group. From now on we will always assume that $D$ is regular. 

To maintain symmetry between the three ramification points of $\beta$, it is sometimes useful to extend $D$ to a tripartite triangulation $\beta^{-1}({\mathbb P}^1({\mathbb R}))$ of $C$ by adding a generator $Z$ to $\Gamma(2)$ satisfying $XYZ=1$, so that the corresponding element $z=(xy)^{-1}$ of $G$ represents the monodromy permutation of $\beta$ around $\infty$. We say that $D$ has {\it type} $(l, m, n)$ where $l, m$ and $n$ are the orders of $x, y$ and $z$. Then the Riemann-Hurwitz formula implies that $D$ has genus
\[g=1+\frac{1}{2}\left(1-\frac{1}{l}-\frac{1}{m}-\frac{1}{n}\right)|G|.\]

\section{Counting regular dessins}

The correspondences in Section~\ref{des} provide a purely group-theoretic approach to dessins. In particular, they give a bijection between isomorphism classes of regular dessins $D$ and normal subgroups $N$ of finite index in a free group
\[F=F_2=\langle X, Y\mid - \rangle = \langle X, Y, Z \mid XYZ=1 \rangle\]
of rank $2$, with ${\rm Aut}\,D\cong F/N$. For each $2$-generator finite group $G$, let ${\mathcal R}(G)$ denote the set of (isomorphism classes of) regular dessins $D$ with ${\rm Aut}\,D\cong G$, and let
\[{\mathcal N}(G)=\{N\triangleleft F \mid F/N\cong G\},\]
so that the above bijection restricts to a bijection between the finite sets ${\mathcal N}(G)$ and ${\mathcal R}(G)$, with each $N\in{\mathcal N}(G)$ corresponding to a dessin $D=D(N)\in{\mathcal R}(G)$. (The set corresponding to ${\mathcal N}(G)$, with $F$ a free group of arbitrary rank, was studied by Dunwoody in~\cite{DunRG}.)

There is also a natural bijection between ${\mathcal N}(G)$ and the set of orbits of ${\rm Aut}\,G$ on generating pairs $(x,y)$ for $G$, or equivalently on generating triples $(x,y,z)$ for $G$ satisfying $xyz=1$. This is because epimorphisms $F\to G$ correspond to choices of generating pairs $(x,y)$ for $G$, and two such epimorphisms have the same kernel $N$ if and only if they differ by an automorphism of $G$. The number
\[r=r(G)=|{\mathcal R}(G)|\]
of regular dessins associated with $G$ is therefore equal to the number
\[d_2(G)=\frac{\phi_2(G)}{|{\rm Aut}\,G|}\]
of such orbits, where $\phi_2(G)$ is the number of generating pairs $(x,y)$ for $G$ (since these are permuted semiregularly by ${\rm Aut}\,G$).

Hall~\cite{Hal} developed a very useful method for evaluating functions such as $\phi_2(G)$. Since every pair of elements of $G$ generate a unique subgroup $H$ we have
\[|G|^2=\sum_{H\le G}\phi_2(H).\]
From this, M\"obius inversion in the subgroup lattice of $G$ shows that
\[\phi_2(G)=\sum_{H\le G}\mu_G(H)|H|^2,\]
where $\mu_G$ is the M\"obius function for $G$, defined recursively by
\[\sum_{K\ge H}\mu_G(K)=\delta_{H,G},\]
and $\delta_{H,G}$ denotes the Kronecker delta, equal to $1$ if $H=G$ and $0$ otherwise.
It follows that
\[r(G)=\frac{1}{|{\rm Aut}\,G|}\sum_{H\le G}\mu_G(H)|H|^2.\]

(This method can be extended to count normal subgroups of {\sl any\/} finitely generated group $F$ with a given finite quotient group $G$: see, for instance~\cite{Jon94, JS} where $F$ is a triangle group, and~\cite{Jon95} where $F$ is a surface group.)

\medskip

\noindent{\bf Example 3.1} It is easily seen that $r(C_n)$ is a multiplicative function of $n$, that is, if $m$ and $n$ are coprime then $r(C_{mn})=r(C_m)r(C_n)$. It follows, by considering the case where $n$ is a prime power, that
\[r(C_n)=n\prod_{p|n}\left(1+\frac{1}{p}\right),\]
where $p$ ranges over the distinct primes dividing $n$.

\medskip

\noindent{\bf Example 3.2} If $G$ is the dihedral group $D_n$ of order $2n$, with $n>2$, then there are, up to automorphisms, just three generating triples $(x,y,z)$ for $G$ with $xyz=1$: in each case, two of $x, y$ and $z$ are reflections, and the third is a rotation of order $n$. Thus $r(D_n)=3$. However, if $n=2$ then $G$ is a Klein four-group $V_4=C_2\times C_2$, with $r(D_2)=1$.

\medskip

\noindent{\bf Example 3.3} In~\cite{Hal}, Hall computed the M\"obius function $\mu_G$ for a number of groups $G$, including the simple groups $L_2(p)=PSL_2(p)$ for primes $p\ge 5$. He showed that $r(A_4)=4$, $r(S_4)=9$, $r(A_5)=19$, $r(\hat A_5)=76$ and $r(A_6)=53$, where $\hat A_5$ is the binary icosahedral group, isomorphic to $SL_2(5)$, and that
\[r(L_2(p))=\frac{1}{4}(p+1)(p^2-2p-1)-\epsilon,\]
where $\epsilon = 49, 40, 11$ or $2$ as $p\equiv \pm 1$ and $\pm 1$, or $p\equiv \pm 1$ and $\pm 3$, or $p\equiv \pm 2$ and $\pm 1$, or $p\equiv \pm 2$ and $\pm 3$ mod~$(5)$ and mod~$(8)$. Thus $r(L_2(p))=19$, $57$, $254$ and $495$ for $p=5, 7, 11$ and $13$ respectively. (Note that $L_2(5)\cong A_5$ and $L_2(9)\cong A_6$, so for these groups $r(G)=19$ and $53$ respectively.)

\medskip

\noindent{\bf Example 3.3} Downs~\cite{DowPhD} extended Hall's computation of $\mu_G$ to $L_2(q)$ and $PGL_2(q)$ for all prime powers $q$; see~\cite{DowJLMS} for a proof for $L_2(2^e)$ and a statement of results for $L_2(q)$ where $q$ is odd, and~\cite{DJ} for some combinatorial applications by Downs and the author. The general results are rather complicated, but when $G=L_2(2^e)$ with $e>1$ we find that
\[r(G)=\frac{1}{e}\sum_{f|e}\mu\left(\frac{e}{f}\right)2^f(2^{2f}-2^f-3)\]
where $\mu$ is the classical M\"obius function on $\mathbb N$ (see~\cite{DJ}). Thus for $e=2$ and $3$ we have $r(G)=19$ (as expected, since $L_2(4)\cong A_5$) and $142$. Downs and the author~\cite{DJ13}  have recently obtained the similar formula
\[r(G)=\frac{1}{e}\sum_{f|e}\mu\left(\frac{e}{f}\right)2^f(2^{4f}-2^{3f}-9)\]
where $G$ is the simple Suzuki group $Sz(2^e)$~\cite{ATLAS, Suz, Wil} for some odd $e>1$.

\medskip

For any finite group $G$ we have $\phi_2(G)\le|G|^2$, so
\[r(G)\le\frac{|G|^2}{|{\rm Aut}\,G|}.\]
In particular, if $G$ has trivial centre, so that ${\rm Inn}\,G\cong G$, then
\[r(G)\le\frac{|G|}{|{\rm Out}\,G|}\]
where ${\rm Out}\,G:={\rm Aut}\,G/{\rm Inn}\,G$ is the outer automorphism group of $G$. As a consequence of the classification of finite simple groups, it is known that any such group can be generated by two elements. In fact, results of Dixon~\cite{Dix}, of Kantor and Lubotzky~\cite{KL}, and of Liebeck and Shalev~\cite{LS} show that if $G$ is a non-abelian finite simple group, then a randomly-chosen pair of elements generate $G$ with probability approaching $1$ as $|G|\to\infty$, so for such groups this upper bound is asymptotically sharp, that is,
\[r(G)\sim\frac{|G|}{|{\rm Out}\,G|}\quad{\rm as}\quad |G|\to\infty.\]
The values of $|G|$ and $|{\rm Out}\,G|$ given in~\cite{ATLAS, Wil} show that for each of the infinite families of non-abelian finite simple groups, $|{\rm Out}\,G|$ grows much more slowly than $|G|$ (indeed, in many cases it is bounded), so that $r(G)$ grows almost as quickly as $|G|$. For instance, $r(A_n)\sim n!/4$ as $n\to\infty$. (See~\cite{Dix05, MT} for more precise results concerning generating pairs for $A_n$.)

One can illustrate how close this asymptotic estimate is as follows. For any finite group $G$ we have
\[\phi_2(G)=|G^2\setminus\bigcup_MM^2|\ge |G|^2-\sum_M|M|^2,\]
where $M$ ranges over the maximal subgroups of $G$. If $G$ is perfect then each such $M$ has $|G:M|$ conjugates, so
\[\phi_2(G)\ge|G|^2\left(1-\sum_{i=1}^k\frac{1}{|G:M_i|}\right),\]
where $M_i$ ranges over a set of representatives of the $k$ conjugacy classes of maximal subgroups of $G$. It follows that if $G$ also has trivial centre (and in particular if $G$ is a non-abelian finite simple group) then
\[1\ge\frac{r(G)}{|G|/|{\rm Out}\,G|}\ge 1-\sum_{i=1}^k\frac{1}{|G:M_i|}.\]
When $G$ is large, the sum on the right is typically very small, so that
\[r(G)\approx\frac{|G|}{|{\rm Out}\,G|}.\]
Thus, for the simple groups $G=L_2(q)$, $Sz(q)$, $R(q)$ and $U_3(q)$ (see~\cite{ATLAS}) we have
\[\sum_{i=1}^k\frac{1}{|G:M_i|}=\frac{1}{q^m}+O\left(\frac{1}{q^{m+1}}\right)
\quad{\rm as}\quad q\to\infty,\]
where $m=1, 2, 3, 3$ respectively, with the sum dominated by the term corresponding to the doubly transitive permutation representation of degree $q^m+1$. The situation is similar for the $26$ sporadic simple groups: for instance, if $G$ is O'Nan's group $O$'$N$~\cite{ONan} (the 13th largest of these groups), with $k=13$ conjugacy classes of maximal subgroups~\cite{ATLAS, Wil85, Yos}, we have
\[\frac{|G|}{|{\rm Out}\,G|}=\frac{460815505920}{2}=230407752960\]
and
\[\sum_{i=1}^k\frac{1}{|G:M_i|}=0.00001726863378\ldots,\]
so
\[230403774132\le r(G)\le 230407752960.\]
If $G$ is the Monster, the largest sporadic simple group, then $|{\rm Out}\,G|=1$ and so
\[r(G)\approx |G|\approx 8.08017\times 10^{53}.\]

\section{Classifications}

For some groups $G$, the dessins in ${\mathcal R}(G)$ have been classified.

\medskip

\noindent{\bf Example 4.1} If $G=C_n$ then three of the dessins in ${\mathcal R}(G)$ have genus $0$, obtained from epimorphisms $F\to G$ sending one of $X, Y$ and $Z$ to the identity; these correspond to Bely\u\i\/ functions $C={\mathbb P}^1({\mathbb C})\to{\mathbb P}^1({\mathbb C})$ given by $\beta:t\mapsto t^n$, $1/(1-\beta)$ and $(\beta-1)/\beta$, resulting in dessins of types $(n,1,n)$, $(n,n,1)$ and $(1,n,n)$ respectively; they have vertices, of three different colours, at $0$, at $\infty$, and at the $n$th roots of $1$. The remaining dessins depend on the arithmetic nature of $n$. In the simplest case, when $n$ is an odd prime $p$, they all have type $(p,p,p)$, they lie on the Lefschetz curves of genus $(p-1)/2$, with affine models $u^p=v^m(v-1)$ for $m=1,\ldots, p-2$, and they have Bely\u\i\/ functions $\beta:(u,v)\mapsto v$.

\medskip

\noindent{\bf Example 4.2} If $G=D_n$ with $n>2$ then the three dessins in ${\mathcal R}(G)$ all have genus $0$, and their types are permutations of $(2,2,n)$. They have vertices in ${\mathbb P}^1({\mathbb C})$, of three different colours, at $0$ and $\infty$, at the $n$th roots of $1$, and at the $n$th roots of $-1$.

\medskip

\noindent{\bf Example 4.3} For $G=A_5$, the $r(G)=19$ dessins in ${\mathcal R}(G)$ have been described (as oriented hypermaps) by Breda and the author in~\cite{BJ} (see~\cite{NN} for representative generating pairs). There are six of genus $0$, each with type a permutation of $(2,3,5)$; there are three each of genera $4$ and $5$, with permutations of $(2,5,5)$ and $(3,3,5)$ as types; there are six of genus $9$, with permutations of $(3,5,5)$ as types (three having the generators of order $5$ conjugate in $G$, and three with them not conjugate), and there is one of genus $13$ and type $(5,5,5)$.

\medskip

\noindent{\bf Example 4.4} Let $G$ be the non-abelian group of order $pq$, where $p$ and $q$ are primes with $p\equiv 1$ mod~$(q)$. This is a semidirect product of a normal subgroup $P\cong C_p$ by a complement $Q\cong C_q$, and it has automorphism group ${\rm Aut}\,G\cong AGL_1(p)$ of order $p(p-1)$. Hall's method gives $r(G)=d_2(G)=q^2-1$. The dessins in ${\mathcal R}(G)$ have been investigated by Streit and Wolfart in~\cite{SW} and (together with the present author) in~\cite{JSW}. There are $3(q-1)$ of genus $(pq-3p+2)/2$, whose types are permutations of $(q,q,p)$, and $(q-1)(q-2)$ of genus $(p-1)(q-2)/2$ and type $(q,q,q)$. These dessins are defined over the cyclotomic field ${\mathbb Q}(\zeta_q)$, where $\zeta_q:=\exp(2\pi i/q)$, with the orbits of the absolute Galois group $\Gamma={\rm Gal}\,\overline{\mathbb Q}/{\mathbb Q}$ on them all having length $q-1$.

\section{Universal covers}

For any $2$-generator finite group $G$, let
\[K(G)=\negthinspace\negthinspace\bigcap_{N\in{\mathcal N}(G)}
\negthinspace\negthinspace\negthinspace N,\]
the intersection of all $N\triangleleft\, G$ with $F/N\cong G$. As the intersection of finitely many normal subgroups of finite index in $F$, $K(G)$ is also a normal subgroup of finite index in $F$. It corresponds to a regular dessin
\[U(G)=D(K(G))=\negthinspace\negthinspace\bigvee_{D\in{\mathcal R}(G)}
\negthinspace\negthinspace\negthinspace D\]
called the {\it universal cover} of $G$, the smallest dessin covering each $D\in{\mathcal R}(G)$, with automorphism group
\[\overline G:={\rm Aut}\,U(G)\cong F/K(G).\]

The dessin $U(G)$ has some special properties. For instance, by its construction, $K(G)$ is the only normal subgroup of $F$ with quotient isomorphic to $\overline G$. This uniquenessss implies that $U(G)$ is invariant under the action of the absolute Galois group $\Gamma={\rm Gal}\,\overline{\mathbb Q}/{\mathbb Q}$, and is therefore defined over $\mathbb Q$.

\medskip

\noindent{\bf Example 5.1} If $G=C_n$ then
\[K(G)=F'F^n,\]
the group generated by the commutators and $n$th powers of elements of $F$, so
\[\overline G=F/F'F^n\cong C_n\times C_n,\]
and $U(G)$ is the $n$th degree Fermat dessin, on the Fermat curve
\[x_0^n+x_1^n+x_2^n=0\]
in ${\mathbb P}^2({\mathbb C})$. Here $\overline G$ acts by multiplying homogeneous coordinates $x_i$ by $n$th roots of $1$. (In fact, we obtain the same universal cover $U(G)$ and group $\overline G$ whenever $G$ is a $2$-generator abelian group of exponent $n$.)

\medskip

\noindent{\bf Example 5.2} If $G=D_n$ with $n>2$ then $\overline G$ is an extension of a normal subgroup $C_n^3$ by a Klein four-group $D_2$, each of its three involutions centralising one of the three direct factors and inverting the other two. (Herrlich~\cite[Prop.~4.5]{Her} has also constructed this covering group in the context of origamis.) The dessin $U(G)$ is an $n^3$-sheeted regular covering of $U(D_2)$, branched over its vertices $\pm 1, \pm i, 0$ and $\infty$. It has type $(e,e,e)$ where $G$ has exponent $e={\rm lcm}\,(n,2)$, and its genus is $2n^3-6n^2+1$ or $2n^3-3n^2+1$ as $n$ is even or odd.

\medskip

\noindent{\bf Example 5.3} If $G$ is the non-abelian group of order $pq$ in Example~4.4, then $\overline G$ is a semidirect product of a normal subgroup $\tilde P\cong C_p^{q^2-1}$ by a complement $\overline Q\cong C_q^2$. As an ${\mathbb F}_p\overline Q$-module, $\tilde P$ is the quotient of the regular module ${\mathbb F}_p\overline Q$ by the principal submodule, i.e.~the direct sum of all the ($1$-dimensional) non-principal irreducible ${\mathbb F}_p\overline Q$-modules. Factoring out all but one of these $q^2-1$ submodules gives a quotient of $\overline G$ isomorphic to $G\times C_q$, and factoring out the second direct factor of this yields $G$. The $q^2-1$ equivalence classes of epimorphisms $F\to G$ all arise in this way. The dessin $U(G)$ has type $(pq,pq,pq)$ and genus
\[g=1+\frac{1}{2}p^{q^2-2}q(pq-3).\]
It is a regular covering of the Fermat dessin $U(Q)$ of degree $q$, branched over its $3q$ black, white and red vertices.

\medskip

The group $\overline G=F/K(G)$ is naturally embedded in the cartesian product
\[\prod_{N\in{\mathcal N}(G)}\negthinspace\negthinspace\negthinspace F/N\cong G^r\]
of $r=r(G)$ copies of $G$ as the subgroup of $G^r$ generated by $(x_i)$ and $(y_i)$ where $x_i$ and $y_i$ are the monodromy generators of $G$ corresponding to the $i$th dessin in some numbering of ${\mathcal R}(G)$. It follows that $\overline G$ satisfies all the identical relations satisfied by $G$. For instance, if $G$ is nilpotent of class $c$, is solvable of derived length $d$, or has exponent $e$, then $\overline G$ also has the same property. Since $\overline G$ is also a $2$-generator group, this might suggest that it is not much larger than $G$. The examples in the next section show that this is not always so.

\section{Non-abelian finite simple groups}

The examples $G=C_n$ and $D_n$ show that when $\overline G$ is embedded in $G^r$, it may be a proper subgroup, even though it projects onto each of the $r$ direct factors. However, if $G$ is a non-abelian finite simple group then $\overline G=G^r$. (This result and its proof are probably well-known to most group-theorists, but they are included here for the benefit of non-specialists.) For background on the finite simple groups, see~\cite{ATLAS, Wil}.

\begin{lemma}\label{max}
If $H=H_1\times\cdots\times H_r$ where $H_1,\ldots, H_r$ are non-abelian simple groups, then the maximal normal subgroups of $H$ are those of the form $\prod_{i\ne j}H_i$ where $j=1,\ldots, r$.
\end{lemma}

\noindent{\sl Proof.} Clearly these are maximal normal subgroups of $H$. Conversely, let $N$ be a maximal normal subgroup, so $H_j\not\le N$ for some $j=1,\ldots, r$. Then $H_j\cap N$ is a proper normal subgroup of $H_j$, so $H_j\cap N=1$ by the simplicity of $H_j$. Since $H_j$ and $N$ are both normal in $H$, their commutator $[H_j,N]$ is contained in $H_j\cap N$. Thus $[H_j,N]=1$, so $N$ is contained in the centraliser $C_H(H_j)$ of $H_j$ in $H$, which is $\prod_{i\ne j}H_i$. Hence $N=\prod_{i\ne j}H_i$ by the maximality of $N$.\hfill$\square$

\medskip

Indeed, it follows that the only normal subgroups of $H$ have the form $\prod_{i\in I}H_i$ where $I\subseteq\{1,\ldots, r\}$. Note that finiteness is not assumed in this proof.

\begin{cor}
Let $N_1,\ldots, N_r$ be distinct normal subgroups of a group $\Phi$, with each $H_i:=\Phi/N_i$ non-abelian and simple. If $K=N_1\cap\cdots\cap N_r$ then
\[\Phi/K\cong H_1\times\cdots\times H_r.\]
\end{cor}

\noindent{\sl Proof.} We use induction on $r$. The case $r=1$ is trivial, so let $r>1$, assume the result for intersections of $r-1$ normal subgroups, and define $L=N_1\cap\cdots\cap N_{r-1}$. If $L\le N_r$ then $N_r/L$ is a maximal normal subgroup of the group $\Phi/L\cong H_1\times\cdots\times H_{r-1}$, so Lemma~\ref{max} implies that $N_r=N_i$ for some $i<r$, against our hypothesis. Thus $L\not\le N_r$, so $LN_r$ is a normal subgroup of $\Phi$ properly containing $N_r$, and hence $LN_r=\Phi$ by the maximality of $N_r$. Since $L\cap N_r=K$ we therefore have
\[\Phi/K=(N_r/K)\times(L/K)\cong(\Phi/L)\times(\Phi/N_r)
\cong(H_1\times\cdots\times H_{r-1})\times H_r,\]
as required. \hfill$\square$

\begin{cor}
If $G$ is a non-abelian finite simple group then $\overline G\cong G^r$ where $r=r(G)$.
\end{cor}

\noindent{\sl Proof.} In Corollary~6.2 take $\Phi=F$ and $\{N_1,\ldots, N_r\}={\mathcal N}(G)$, so that $K=K(G)$ and each $H_i\cong G$; then $\overline G$ can be identified with $G^r$, as claimed.\hfill$\square$

\medskip

Guralnick and Kantor~\cite{GK} have shown that if $G$ is a non-abelian finite simple group then every non-identity element is a member of a generating pair. It follows that if such a group $G$ has exponent $e$ then $U(G)$ has type $(e,e,e)$, so by the Riemann-Hurwitz formula~\cite[Remark~1.2.21]{LZ} it has genus
\[g=1+\frac{e-3}{2e}|G|^r.\]
In this case, the numerical estimates given earlier for $r(G)$ show that $\overline G$, and hence $g$, can be rather large.

\medskip

\noindent{\bf Example~6.1} In the case of the smallest non-abelian finite simple group, namely $G=A_5$, we have $r(G)=19$, so $\overline G=G^{19}$, of order
\[60^{19}=609359740010496\times 10^{17}\approx 6.1\times 10^{31}.\]
Since $G$ has exponent $30$, it follows that $U(G)$ has type $(30,30,30)$ and genus
\[1+\frac{9}{20}\times 60^{19}=274218830047232000000000000000001\approx2.742\times 10^{31}.\]

\medskip

\noindent{\bf Example~6.2} For the second smallest non-abelian finite simple group, $G=L_2(7)$, we have $r(G)=57$, so $\overline G=G^{57}$, of order $168^{57}\approx 7.035\times 10^{126}$. The universal cover $U(G)$ has type $(84,84,84)$ and genus
\[1+\frac{27}{56}\times 168^{57}\approx 3.392\times 10^{126}.\]

\medskip

\noindent{\bf Example~6.3} The Monster simple group $G$ has order
\[|G|=2^{46}.3^{20}.5^9.7^6.11^2.13^3.17.19.23.29.31.41.47.59.71\]
\[=808017424794512875886459904961710757005754368000000000\]
\[\approx 8.080\times 10^{53}\]
and exponent
\[e=2^5.3^3.5^2.7.11.13.17.19.23.29.31.41.47.59.71\]
\[=1165654792878376600800\approx 1.166\times 10^{21}. \]
Since $|{\rm Out}\,G|=1$ we have $r(G)\approx|G|$, so
\[|\overline G|=|G|^r\approx|G|^{|G|}\approx(8.080\times 10^{53})^{8.080\times 10^{53}}\approx 10^{10^{55.639}}.\]
The universal cover $U(G)$ has type $(e,e,e)$ and genus approximately $|\overline G|/2$.

\medskip

\noindent{\bf Example~6.4} For comparison, we consider a quasisimple example. Let $G$ be the binary icosahedral group $\hat A_5\cong SL_2(5)$, with ${\rm Aut}\,G\cong S_5\cong PGL_2(5)\cong{\rm Aut}\,A_5$. Then $r(G)=76$, with each generating triple for the central quotient $G/Z(G)\cong A_5$ lifting to four for $G$: the members of each such quadruple are related to each other by multiplying two members of a triple by the central involution in $G$. Each generating triple for $G$ arises in this way, so the $r(A_5)=19$ orbits of $S_5$ on generating pairs for $A_5$ yield $r(G)=4\times 19=76$ orbits on those for $G$. Each $N\in{\mathcal N}(A_5)$ contains a unique normal subgroup $N^*=N\cap F'F^2$ of $F$ with $F/N^*\cong G\times V_4$; this quotient has an elementary abelian centre $E$ of order $8$, and the four maximal subgroups of $E$ complementing $Z(G)$ lift to the four subgroups of $F$ in ${\mathcal N}(G)$ corresponding to $N$. Using this, one can show that $K(G)$ is a subgroup of index $2$ in $K(A_5)\cap F'F^2$, and that $\overline G=C\times V_4$ where $C$ is the central product of $19$ copies of $G$ (the quotient of $G^{19}$ formed by identifying the centres of the direct factors). The dessin $U(G)$, an $8$-sheeted branched covering of $U(A_5)$, has type $(60,60,60)$ and genus
\[231556701203988480000000000000001\approx 2.316\times 10^{32}.\]

\section{Operations on dessins and $T_2$-systems}

The automorphism group of the free group $F$ permutes the subgroups of $F$. Since inner automorphisms leave invariant each conjugacy class of subgroups, this induces an action of the outer automorphism group
\[\Omega:={\rm Out}\,F={\rm Aut}\,F/{\rm Inn}\,F\]
of $F$ on isomorphism classes of dessins. In~\cite{Jam}, James interpreted $\Omega$ as the group of all operations on oriented hypermaps (which, in the finite case, are equivalent to dessins); see~\cite{JSin, JSW} for operations on various other categories of maps, hypermaps and dessins. Here we will consider, for each finite group $G$, the isomorphic actions of $\Omega$ on ${\mathcal N}(G)$ and on ${\mathcal R}(G)$.

For any integer $n\ge 1$, the automorphism group of the free group $F_n$ of rank $n$ is generated by the elementary Nielsen transformations: permuting the free generators, inverting one of them, and multiplying one of them by another~\cite[Theorem~3.2]{MKS}. When $n=2$, with $F_2=F$, one can identify $\Omega={\rm Out}\,F$ with $GL_2({\mathbb Z})$ through its faithful induced action on the abelianisation $F^{\rm ab}=F/F'\cong{\mathbb Z}^2$ of $F$~\cite[Ch.~I, Prop.~4.5]{LSch}. If we take the images of $X$ and $Y$ in $F^{\rm ab}$ as a basis, the elementary Nielsen transformations are represented by the matrices
\[E_i=
\left(\,\begin{matrix}0&1\cr 1&0\cr \end{matrix}\,\right),\quad
\left(\,\begin{matrix}-1&0\cr 0&1\cr \end{matrix}\,\right),\quad
\left(\,\begin{matrix}1&0\cr 0&-1\cr \end{matrix}\,\right),\quad
\left(\,\begin{matrix}1&1\cr 0&1\cr \end{matrix}\,\right)\quad{\rm and}\quad
\left(\,\begin{matrix}1&0\cr 1&1\cr \end{matrix}\,\right)\]
for $i=1,\ldots, 5$. This is not a minimal set of generators: for instance, $E_1$ and $E_4$ suffice to generate $\Omega$, since $\langle E_4, E_5\rangle=SL_2({\mathbb Z})$ and $E_5=E_1E_4E_1$.)

The element $E_1$ of $\Omega$ is represented by the automorphism of $F$ transposing $X$ and $Y$. This corresponds to the operation of interchanging the colours of the black and white vertices of a dessin, and is equivalent to replacing the Bely\u\i\/ function $\beta$ with $1-\beta$, or to transposing the monodromy generators $x$ and $y$ of $G$. This sends $Z=(XY)^{-1}$ to $Z^X$, so it sends a dessin of type $(l,m,n)$ to one of type $(m,l,n)$ on the same surface. The element
\[S:=E_4^{-1}E_5=\left(\,\begin{matrix}0&-1\cr 1&1\cr \end{matrix}\,\right)\]
has order $6$; this lifts only to automorphisms of infinite order~\cite{JP}, but
\[T:=S^2=(E_4^{-1}E_5)^2=\left(\,\begin{matrix}-1&-1\cr 1&0\cr \end{matrix}\,\right)\]
lifts to the automorphism $X\mapsto Z\mapsto Y\mapsto X$ which permutes the vertex colours black, white and red in a $3$-cycle. The elements $T$ and $E_1$ generate a subgroup $\Sigma$ of $\Omega$ isomorphic to $S_3$, described by Mach\`\i\/ in~\cite{Mac}, inducing all permutations of the three colours and preserving the underlying surface. The element $S^3$ is the central involution $-I$ of $\Omega$, induced by the automorphism inverting $X$ and $Y$ and reversing the orientation of each dessin. (Factoring out $\langle -I\rangle$ gives the central quotient $\Omega/Z(\Omega)\cong PGL_2({\mathbb Z})$, the group of operations on all hypermaps, ignoring orientation~\cite{Jam, JSin}.) The subgroup $\Omega_1:=\langle E_1, S\rangle=\Sigma\times\langle -I\rangle\cong S_3\times C_2\cong D_6$ of $\Omega$ preserves the genus of each dessin, but its elements may reverse the orientation and permute the vertex colours.

The element
\[U:=SE_4^{-1}=\left(\,\begin{matrix}0&-1\cr 1&0\cr \end{matrix}\,\right)\]
corresponds to the automorphism $X\mapsto Y^{-1}$, $Y\mapsto X$ of order $4$. This and its inverse transpose the colours of the black and white vertices, while reversing the cyclic permutations of edges around the vertices of one colour. Unlike the operations in $\Omega_1$, these can (and usually do) change the genus of a dessin, even though the monodromy and automorphism groups are preserved. The elements $U$ and $E_1$ generate a subgroup $\Omega_2\cong D_4$; in addition to $U^2=-I$, this contains the elements
\[E_2=UE_1=\left(\,\begin{matrix}-1&0\cr 0&1\cr \end{matrix}\,\right)
\quad{\rm and}\quad
E_3=E_1U=\left(\,\begin{matrix}1&0\cr 0&-1\cr \end{matrix}\,\right)\]
represented by the automorphisms of $F$ inverting $X$ or $Y$ respectively. The operations corresponding to $E_2$ and $E_3$ reverse the cyclic permutations of edges around the black or white vertices of a dessin; they are sometimes called Petrie operations, since they transpose faces and Petrie polygons (closed zig-zag paths). These are special cases of operations $H_{i,j}:x\mapsto x^i, y\mapsto y^j$ on dessins considered by Streit, Wolfart and the author in~\cite{JSW}, and also of certain operations on regular maps introduced earlier by Wilson in~\cite{WilOp}.

The subgroups $\Omega_1$ and $\Omega_2$ generate $\Omega$. In fact, $\Omega$ is the free product
\[\Omega=\Omega_1*_{\Omega_0}\Omega_2\cong D_6*_{D_2}D_4\]
of these two subgroups, amalgamating their common subgroup
\[\Omega_0:=\langle E_1, -I\rangle\cong D_2\cong V_4.\]
Thus $\Omega$ is generated by operations of finite order; these have been classified by Pinto and the author in~\cite{JP}.

The orbits of $\Omega$ on ${\mathcal R}(G)$, or equivalently on ${\mathcal N}(G)$, correspond to the $T_2$-systems in $G$, that is, the orbits of ${\rm Aut}\,F\times{\rm Aut}\,G$ acting by composition on epimorphisms $F\to G$ and thus on generating pairs $(x,y)$ for $G$. (These are called $T$-systems when the common size of the generating sets and the rank of the free group is arbitrary.) Let $\nu=\nu(G)$ denote the number of orbits of $\Omega$ in each of these actions. Then ${\mathcal N}(G)$ and ${\mathcal R}(G)$ each split into $\nu(G)$ orbits ${\mathcal N}_i(G)$ and ${\mathcal R}_i(G)$ of length $r_i$ for $i=1,\ldots, \nu$, with $r_1+\cdots+r_{\nu}=r(G)$. Each orbit determines a characteristic subgroup
\[K_i(G)=\negthinspace\negthinspace\bigcap_{N\in{\mathcal N}_i(G)}
\negthinspace\negthinspace\negthinspace N\]
of $F$, corresponding to an $\Omega$-invariant regular dessin
\[U_i(G)=D(K_i(G))=\negthinspace\negthinspace\bigvee_{D\in{\mathcal R}_i(G)}
\negthinspace\negthinspace\negthinspace D\]
with automorphism group
${\rm Aut}\,U_i(G)\cong F/K_i(G)$.
We have
\[K(G)=\bigcap_{i=1}^{\nu(G)}K_i(G)
\quad{\rm and}\quad
U(G)=\bigvee_{i=1}^{\nu(G)}U_i(G).\]

\medskip

\noindent{\bf Example 7.1} It is known~\cite{DunTS, NN} that if $G$ is abelian then $\nu(G)=1$, so all dessins in ${\mathcal R}(G)$ are equivalent under $\Omega$. The same applies if $G$ is a dihedral group. On the other hand, Dunwoody~\cite{DunTS} has shown that for each prime $p$ there are $2$-generator $p$-groups $G$ of nilpotence class $2$ with $\nu(G)$ arbitrarily large.

\medskip

An observation of Graham Higman, based on a theorem of Nielsen~\cite{Nie} (see also~\cite[Theorem~3.9]{MKS}), states that if $G=\langle x, y\rangle$ then the union of the orbits of ${\rm Aut}\,G$ containing $c:=[x,y]$ and $c^{-1}=[y,x]$ is invariant under $\Omega$, and hence so is the order of $c$. In many cases this can be used to give lower bounds for $\nu(G)$.

\medskip

\noindent{\bf Example 7.2}  Neumann and Neumann~\cite{NN} showed that if $G=A_5$ then the $19$ elements of ${\mathcal N}(G)$ form two orbits of $\Omega$, of lengths $9$ and $10$, so the same happens to the corresponding dessins. In fact, those of type a permutation of $(2,5,5)$, $(3,3,5)$ or $(3,5,5)^{-}$ form an orbit ${\mathcal R}_1(G)$ of length $r_1=3+3+3=9$, while those of type a permutation of $(2,3,5)$, $(3,5,5)^+$ or $(5,5,5)$ form an orbit ${\mathcal R}_2(G)$ of length $r_2=6+3+1=10$; here the superscript $+$ or $-$ indicates that the generators of order $5$ are or are not conjugate in $G$. These orbits are distinguished by the property that, for each dessin, the commutator $[x,y]$ has order $3$ or $5$ respectively. They correspond to characteristic subgroups $K_i(G)$ of $F$ with quotients $G^{r_i}$, such that $K(G)=K_1(G)\cap K_2(G)$, and to $\Omega$-invariant regular dessins $U_i(G)$ of type $(30,30,30)$ and genera
\[1+\frac{9}{20}\times 60^{r_i}=4934963200000001
\quad{\rm and}\quad 296097792000000001,\]
with $U(G)=U_1(G)\vee U_2(G)$.

\medskip

Garion and Shalev~\cite{GS} have shown that if $G$ is a non-abelian finite simple group, then $\nu(G)\to\infty$ as $|G|\to\infty$, as conjectured by Guralnick and Pak~\cite{GP}.

\medskip

\noindent{\bf Example 7.3} Let $G=L_2(p)$ for some prime $p\ge 5$. Then it follows from Dickson's description~\cite[Ch.~XII]{Dic} of the subgroups of $G$ that the elements
\[x=\pm\left(\,\begin{matrix}1&a\cr 0&1\cr \end{matrix}\,\right)
\quad{\rm and}\quad
y=\pm\left(\,\begin{matrix}1&0\cr 1&1\cr \end{matrix}\,\right)\]
generate $G$ provided $a\ne 0$. This generating pair corresponds to a regular dessin $D_a\in{\mathcal R}(G)$ of type $(p, p, k)$, where $k$ is the order in $G$ of an element with trace $\pm(a+2)$. These generators have commutator
\[c=[x,y]=\pm\left(\,\begin{matrix}a^2+a+1&a^2\cr -a&1-a\cr \end{matrix}\,\right),\]
so $c$ and $c^{-1}$ have trace $\pm(a^2+2)$. Higman's commutator criterion therefore implies that if $\pm(a^2+2)\ne\pm(b^2+2)$ then $D_a$ and $D_b$ lie in different orbits of $\Omega$, so $\nu(G)\ge(p-1)/4$. 

\medskip

An interesting invariant of a regular dessin $D$ is the permutation group induced by $\Omega$ on the orbit containing $D$.

\medskip

\noindent{\bf Example 7.4} According to Neumann and Neumann~\cite{NN}, $\Omega$ acts on ${\mathcal R}(A_5)$ as $S_9\times S_{10}$. This is the largest possible group with orbits of lengths $9$ and $10$, so neither of the two orbits ${\mathcal R}_i(A_5)$ admits any non-trivial $\Omega$-invariant structure.

\medskip

\noindent{\bf Example 7.5} If $G=C_n$ then $\Omega$ acts on $F/K(G)=F/F'F^n\cong{\mathbb Z}_n^2$ as the subgroup $GL_2^{\pm 1}({\mathbb Z}_n)$ of $GL_2({\mathbb Z}_n)$ consisting of the matrices of determinant $\pm 1$. (To see this, simply reduce the matrices $E_i$ modulo~$(n)$, and note that the images of $E_4$ and $E_5$ generate $SL_2({\mathbb Z}_n)$.) It follows that $\Omega$ acts on ${\mathcal R}(C_n)$ as $GL_2^{\pm 1}({\mathbb Z}_n)$ acts on the subgroups of ${\mathbb Z}^2$ with quotient ${\mathbb Z}_n\cong C_n$. For instance, if $n$ is a prime $p$ then ${\mathcal R}(G)$ can be identified with the projective line ${\mathbb P}^1({\mathbb F}_p)$, and $\Omega$ acts on it as $L_2(p)$ if $p=2$ or $p\equiv 1$ mod~$(4)$ (so that $-1$ is a square in ${\mathbb F}_p$), but as $PGL_2(p)$ otherwise. In the latter case, the triple transitivity of $PGL_2(p)$ ensures that there are no $\Omega$-invariant binary or ternary relations on ${\mathcal R}(C_n)$; however, the quaternary relation induced by evaluating the cross-ratio is invariant in both cases.

\medskip

\noindent{\bf Example 7.6} Let $G$ be the non-abelian group of order $pq$ in Examples~4.4 and 5.3. Since $xyz=1$, the images $(\overline x,\overline y,\overline z)$ in $(G/P)^3\cong C_q^3\cong{\mathbb F}_q^3$ of the generating triples $(x,y,z)$ for $G$ are the $q^2-1$ non-identity elements of a subgroup $S\cong F/F'F^q\cong C_q^2\cong{\mathbb F}_q^2$. These are fixed by ${\rm Aut}\,G$, whereas $\Omega$ acts on them as $GL_2^{\pm 1}(q)$ (see Example~7.5), so $\Omega$ acts transitively on ${\mathcal R}(G)$.

If $q=2$ then $\Omega$ acts $3$-transitively on ${\mathcal R}(G)$ as $S_3$, but if $q>2$ then its action is imprimitive, with the non-zero elements of the $1$-dimensional linear subspaces of $S$ forming $q+1$ blocks of size $q-1$; these are permuted transitively by $\Omega$, acting as $L_2(q)$ or $PGL_2(q)$ on ${\mathbb P}^1({\mathbb F}_q)$ as $q\equiv 1$ or $-1$ mod~$(4)$ (again, see Example~7.5). These blocks are the orbits of $\Gamma={\rm Gal}\,\overline{\mathbb Q}/{\mathbb Q}$, which acts on ${\mathcal R}(G)$ as ${\rm Gal}\,{\mathbb Q}(\zeta_q)/{\mathbb Q}\cong{\mathbb F}_q^*$, inducing the group ${\mathbb F}_q^*.I\cong C_{q-1}$ of scalar matrices on $S$. Thus the actions of $\Omega$ and $\Gamma$ on ${\mathcal R}(G)$ commute, generating a group $\langle GL_1^{\pm 1}(q), {\mathbb F}_q^*.I\rangle$ which has index $2$ or $1$ in $GL_2(q)$ as $q\equiv 1$ or $-1$ mod~$(4)$.

\medskip

The actions of $\Gamma$ and $\Omega$ on each of the sets ${\mathcal R}(G)$ considered above commute with each other (in some cases, though not Example~7.6, because one group acts trivially). In fact, the action of $\Gamma$ always commutes with that of the Mach\`\i\/ subgroup $\Sigma$ of $\Omega$, which simply permutes the three vertex colours, but the following example shows that on some sets ${\mathcal R}(G)$ it does not commute with $\Omega$.

\medskip

\noindent{\bf Example 7.7}  There are three non-regular dessins of genus $0$ and type $(6,2,6)$ consisting of plane trees with black and white vertices of valencies $3, 2, 1$ and $2, 2, 1, 1$: one is reflexible and the others form a chiral pair. Defined over the splitting field of the polynomial $25t^3-12t^2-24t-16$, they form an orbit of $\Gamma$, which acts on them as the Galois group $S_3$ of this polynomial (see~\cite[Example~4.58]{GG} or~\cite[\S2.2.2.3]{LZ}, for instance). They each have monodromy group $G\cong S_6$, so their minimal regular covers form an orbit of $\Gamma$ consisting of three dessins of type $(6,2,6)$ and genus $61$ in ${\mathcal R}(S_6)$; again, one is reflexible and the others form a chiral pair. There is an element of $\Gamma$ permuting them in a $3$-cycle, whereas the element $-I$ of $\Omega$ induces a transposition on them, so the actions of $\Gamma$ and $\Omega$ on ${\mathcal R}(S_6)$ do not commute. (In fact, $\langle\Gamma,\Omega\rangle$ acts on its orbit containing these three regular dessins as $S_3\times S_3$ on the $18$ cosets of a subgroup $C_2\times 1$, with $\Gamma$ inducing $S_3\times 1$ and $\Omega$, acting as $\Omega_1=\langle -I\rangle\times\Sigma$, inducing $C_2\times S_3$.)

\section{Concluding remarks}

Dessins can be regarded as unbranched finite coverings of ${\mathbb P}^1({\mathbb C})\setminus\{0, 1, \infty\}$, or equivalently as finite permutation representations of its fundamental group $F=F_2$, so the results obtained here can also be applied to coverings of other spaces with the same fundamental group. These include a bouquet of two circles, with the regular objects the Cayley graphs of $2$-generator groups~\cite{SiW}, and a punctured torus, leading to the theory of origamis: see the survey~\cite{HS} by Herrlich and Schmith\"usen for parallels between origamis and dessins.  In certain other categories, the objects can also be identified with the permutation representations of a particular group: thus one can use the free products $V_4*C_2$ and $C_2*C_2*C_2$ for maps and hypermaps, string Coxeter groups for polytopes, and fundamental groups for covering spaces, etc. The enumerative techniques described here, along with the coverings and operations, can all be applied in similar ways in these contexts: see~\cite{Jon13} for this more general approach.


\begin{thebibliography}{99}

\bibitem{BCIW} P.~Beazley-Cohen, C.~Itzykson and J.~Wolfart, Fuchsian triangle groups and Grothendieck dessins: variations on a theme of Belyi, {\sl~Comm.~Math.~Phys.} 163 (1994), 605--627. 




\bibitem{Bel} G.~V.~Bely\u\i, On Galois extensions of a maximal cyclotomic field, {\sl Izv.~Akad. Nauk SSSR Ser Mat.} 43 (1979), 267--276, 479. 

\bibitem{BJ} A.~J.~Breda d'Azevedo and G.~A.~Jones, Platonic hypermaps, {\sl Beitr\"age Algebra Geom.} 42 (2001), 1--37. 





\bibitem{ATLAS} J.~H.~Conway, R.~T.~Curtis, S.~P.~Norton, R.~A.~Parker and R.~A.~Wilson, {\sl ATLAS of Finite Groups\/}, Clarendon Press, Oxford, 1985. 

\bibitem{Dic} L.~E.~Dickson, {\sl Linear Groups\/}, Dover, New York, 1958. 

\bibitem{Dix} J.~D.~Dixon, The probability of generating the symmetric group, {\sl Math.~Z.} 110 (1969), 199--205. 

\bibitem{Dix05} J.~D.~Dixon, Asymptotics of generating the symmetric and alternating groups, {\sl Electron.~J.~Combin.} 12 (2005), Research paper 5, 5 pp. 

\bibitem{DowPhD} M.~L.~N.~Downs, PhD thesis, University of Southampton, Southampton, UK, 1988. 

\bibitem{DowJLMS} M.~L.~N.~Downs, The M\"obius function of $PSL_2(q)$, with application to the maximal normal subgroups of the modular group, {\sl J.~London Math.~Soc.} 43 (1991), 61--75. 

\bibitem{DJ} M.~L.~N.~Downs and G.~A.~Jones, Enumerating regular objects with a given automorphism group, {\sl Discrete Math.} 64 (1987), 299--302. 

\bibitem{DJ13} M.~L.~N.~Downs and G.~A.~Jones, Enumerating regular objects associated with Suzuki groups, preprint. 

\bibitem{DunRG} M.~J.~Dunwoody, On relation groups, {\sl Math.~Z.} 81 (1963), 180--186. 

\bibitem{DunTS}  M.~J.~Dunwoody, On $T$-systems of groups, {\sl J.~Aust.~Math.~Soc.} 3 (1963), 172--179. 





\bibitem{GS} S.~Garion and A.~Shalev, Commutator maps, measure preservation, and $T$-systems, {\sl Trans.~Amer.~Math.~Soc.} 361 (2009), 4631--4651. 

\bibitem{GG} E.~Girondo and G.~Gonz\'alez-Diez, {\sl Introduction to Compact Riemann Surfaces and Dessins d'Enfants}, London Math.~Soc.~Student Texts 79, Cambbridge University Press, Cambridge, 2011. 

\bibitem{GJ} G.~Gonz\'alez-Diez and A.~Jaikin-Zapirain, The absolute Galois group acts faithfully on regular dessins and on Beauville surfaces, preprint, 2013. 

\bibitem{Gro} A.~Grothendieck, Esquisse d'un Programme, in {\sl Geometric Galois Actions~I, Around Grothendieck's Esquisse d'un Programme} (ed.~P.~Lochak and L.~Schneps) London Math.~Soc.~Lecture Note Ser.~242 (Cambridge University Press, Cambridge, 1997), 5--48. 

\bibitem{GK} R.~M.~Guralnick and W.~M.~Kantor, Probabilistic generation of finite simple groups, {\sl J.~Algebra} 234 (2000), 743--792. 


\bibitem{GP} R.~M.~Guralnick and I.~Pak, On a question of B.~H.~Neumann, {\sl Proc.~Amer. Math.~Soc.} 131 (2002), 2021--2025. 

\bibitem{Hal} P.~Hall, The Eulerian functions of a group, {\sl Q.~J.~Math.} 7 (1936), 134--151. 


\bibitem{Her} F.~Herrlich, Teichm\"uller curves defined by characteristic origamis, in {\sl The Geometry of Riemann Surfaces and Abelian Varieties}, 133--144, {\sl Contemp. Math.} 397, Amer.~Math.~Soc., Providence, RI, 2006. 

\bibitem{HS} F.~Herrlich and G.~Schmith\"usen, Dessins d'enfants and origami curves, in {\sl Handbook of Teichm\"uller Theory, Vol.~II,} 767--809, IRMA Lect.~Math.~Theor.~Phys. 13, Eur.~Math.~Soc., Z\"urich, 2009.



\bibitem{Jam} L.~D.~James, Operations on hypermaps, and outer automorphisms, {\sl European J.~Combin.} 9 (1988), 551--560. 


\bibitem{Jon94} G.~A.~Jones, Ree groups and Riemann surfaces, {\sl J.~Algebra} 165 (1994), 41--62. 

\bibitem{Jon95} G.~A.~Jones, Enumeration of homomorphisms and surface-coverings, {\sl Q.~J.~Math.} 46 (1995), 485--507. 

\bibitem{Jon13} G.~A.~Jones, Combinatorial categories and permutation representations, preprint. 



\bibitem{JP} G.~A.~Jones and D.~Pinto, Hypermap operations of finite order, {\sl Discrete Math.} 310 (2010), 1820--1827. 

\bibitem{JS} G.~A.~Jones and S.~A.~Silver, Suzuki groups and surfaces, {\sl J.~London Math. Soc.} (2) 48 (1993), 117--125. 

\bibitem{JSin} G.~A.~Jones and D.~Singerman, Maps, hypermaps and triangle groups., in {\sl The Grothendieck Theory of Dessins d'Enfants (Luminy, 1993)}, 115--145, London Math.~Soc.~Lecture Note Ser.~200, Cambridge Univ.~Press, Cambridge, 1994. 

\bibitem{JSW} G.~A.~Jones, M.~Streit and J.~Wolfart, Wilson's map operations on regular dessins and cyclotomic fields of definition, {\sl Proc.~Lond.~Math.~Soc.} (3) 100 (2010), 510--532. 


\bibitem{KL} W.~M.~Kantor and A.~Lubotzky, The probability of generating a finite classical group, {\sl Geom.~Dedicata} 36 (1990), 67--87. 



\bibitem{LZ} S.~K.~Lando and A.~K.~Zvonkin, {\sl Graphs on Surfaces and their Applications}, Springer-Verlag, Berlin -- Heidelberg -- New York, 2004. 

\bibitem{LS} M.~W.~Liebeck and A.~Shalev, The probability of generating a finite simple group, {\sl Geom.~Dedicata} 56 (1995), 103--113. 

\bibitem{LSch} R.~C.~Lyndon and P.~E.~Schupp, {\sl Combinatorial Group Theory}, Springer-Verlag, Berlin -- Heidelberg -- New York, 1977. 

\bibitem{Mac} A.~Mach\`\i\/, On the complexity of a hypermap, {\sl Discrete Math.} 42 (1982), 221--226.


\bibitem{MKS} W.~Magnus, A.~Karrass and D.~Solitar, {\sl Combinatorial Group Theory}, Dover, New York, 1976. 

\bibitem{MT} A.~Mar\'oti and M.~C.~Tamburini, Bounds for the probability of generating the symmetric and alternating groups, {\sl Arch.~Math.~(Basel)} 96 (2011), 115--121. 



\bibitem{NN} B.~H.~Neumann and H.~Neumann, Zwei Klassen charakteristischer Untergruppen und ihre Faktorgruppen, {\sl Math.~Nachr.} 4 (1951), 106--125. 

\bibitem{Nie} J.~Nielsen, Die Isomorphismen der allgemeinen unendlichen Gruppe mit zwei Erzeugenden, {\sl Math.~Ann.} 78 (1918), 385--397. 

\bibitem{ONan} M.~E.~O'Nan, Some evidence for the existence of a new simple group, {\sl Proc.~London Math.~Soc.} (3) 32 (1976), 421--479. 




\bibitem{SiW} D.~Singerman and J.~Wolfart, Cayley graphs, Cori hypermaps, and dessins d'enfants, {\sl Ars Math.~Contemp.} 1 (2008), 144--153.

\bibitem{SW} M.~Streit and J.~Wolfart, Characters and Galois invariants of regular dessins, {\sl Rev.~Mat.~Complut.} 13 (2000), 49--81. 

\bibitem{Suz} M.~Suzuki, On a class of doubly transitive groups, {\sl Ann.~of Math.} (2) 75 (1962), 105--145. 


\bibitem{Wil85} R.~A.~Wilson,  The maximal subgroups of the O'Nan group, {\sl J.~Algebra} 97 (1985), 467--473. 

\bibitem{Wil} R.~A.~Wilson, {\sl The Finite Simple Groups}, Graduate Texts in Math. 251, Springer, London, 2009. 

\bibitem{WilOp} S.~E.~Wilson, Operators over regular maps, {\sl Pacific J.~Math.} 81 (1979), 559--568. 

\bibitem{Wol} J.~Wolfart, $ABC$ for polynomials, dessins d'enfants, and uniformization --- a survey, {\sl Elementare und Analytische Zahlentheorie (Tagungsband)}, Proc.~ELAZ-Conference, 24--28 May, 2004 (ed.~W.~Schwarz and J.~Steuding, Steiner, Stuttgart, 2006); http://www.math.uni-frankfurt.de/$\sim$wolfart. 

\bibitem{Yos} S.~Yoshiara, The maximal subgroups of the sporadic simple group of O'Nan, {\sl J.~Fac.~Sci.~Univ.~Tokyo Sect.~IA Math.} 32 (1985), 105--141. 


\end{thebibliography}
\end{document}